\magnification=\magstep1
\hsize=14 truecm
\vsize=21 truecm
\voffset -1.5truecm
\input amstex
\documentstyle{amsppt}

\font\tenscr=rsfs10 
\font\sevenscr=rsfs7 
\font\fivescr=rsfs5 
\skewchar\tenscr='177 \skewchar\sevenscr='177
\skewchar\fivescr='177
\newfam\scrfam \textfont\scrfam=\tenscr \scriptfont\scrfam=\sevenscr
\scriptscriptfont\scrfam=\fivescr
\def\scr{\fam\scrfam}

\def\le{\leqslant}
\def\ge{\geqslant}
\def\leq{\leqslant}

\topmatter
\topinsert
\captionwidth{10 truecm}
\flushpar Acta Math. Appl. Sin.11:4(1995), 389-404.
\endinsert
\title On the Optimality in General Sense for Odd-Block Search\endtitle
\author{Mu-Fa Chen and Dan-Hua Huang}\endauthor
\affil{(Beijing Normal University and Fujian Normal University)}\endaffil
\address{Department of Mathematics, Beijing Normal University, Beijing 100875.
  \newline
  \text{\hskip1em} Department of Mathematics, Fujian Normal University, Fujian Province 350007.}
  \endaddress
\thanks{ This work is partially supported by Ying-Tung Fok Educational Foundation and
Natural Science Foundation of China}\endthanks

\abstract {In his classical article[3](1953), J.Kiefer
introduced the Fibonacci search as a direct optimal  method.  The
optimality was proved under the restriction:  the total number of
tests  is given in advance and fixed.  To avoid this restriction,
some  different  concepts of optimality were  proposed  and  some
corresponding optimal me\-thods were obtained in [1], [2], [5] and [6].
In particular, the  even-block search was treated in [1].
This paper deals with the odd-block search. The main result is Theorem (1.15).}
\endabstract

\endtopmatter

\document

\heading 1. Backgrounds and Main Results  \endheading

In  this  section,  we first review some necessary preliminaries and then
state our main result.
The study of optimal search is usually restricted on the unimodal functions.

\proclaim{(1.1)Definition$^{[3]}$}   A  function  $f$  on  interval
$[0,1]$ is called unimodal if there exists precisely one  maximum
at a point $c_f \in [0,1] $ and the function is strictly monotone
on  $[0,c_f]  $ and $[c_f,1] $. \endproclaim

Let  ${\scr F}  $ denote the set of all  unimodal  functions  on
$[0,1]$.   It is obvious  that  the  unimodal
functions have the following advantage:  Whenever we have had two
tests, we can compare the two results and cancel a part of
the interval. Next, we consider only the following testing methods.

\proclaim{(1.2)Definition$^{[3,8]}$}  A  policy $($or  strategy,  or
sequential search$\,)$ ${\scr P} $ is such a rule:  at the first step, the
rule determines a test point $x_1=x_1( {\scr P} )$  independent
of  $f\in {\scr F} $; at the $n$-th step, the rule determines the $n$-th
test  point  $x_n=x_n( {\scr P},  f)$ according  to  the  first
$(n-1)$ tested points $x_1,\cdots ,x_{n-1} $ and their results
$f(x_1),\cdots ,f(x_{n-1}).$ \endproclaim

An  example  of policies is the {\it Fibonacci search} or
{\it fraction method}  ${\scr F}_n$.  Recall  the
Fibonacci sequence:
$F_0=F_1=1$, $F_n=F_{n-1}+F_{n-2}$, $n\ge 2$.
For  fixed $n\ge 1$,  the  policy ${{\scr F}}_n $ is defined as  follows:
Set $x_1=x_1({{\scr F}}_n)=F_n/F_{n+1}$. Assume that at the  $m$-th
step  $(1\le  m\le  {n-1})$,  we have eliminated a  part  of  the
interval and the remaining interval is $[a_m, b_m]$ with a tested
point $c_m$ inside.  Then, we choose $x_{m+1}=x_{m+1}({{\scr F}}_n,
f)=a_m+b_m-c_m$,  which is just the mirror image of $c_m$ with respect to
the middle point of the interval $[a_m,  b_m]$.  In what follows, we
call the last procedure the {\it symmetry rule}.

\proclaim{Warning}  It  can be happened that for a given  $f$, at some
steps the two tested results are the same. And so after the elimination,
the remaining interval has no tested point inside.  In this  case,
we  have to modify the above ${\scr F}_n$.  Simply regard the remaining
interval as our initial testing interval and apply the same rule.
But in  what  follows we may and will omit this  exceptional  case  for
saving the space. \endproclaim

Let policy ${\scr P}$ act on $f\in {\scr F} $  in $n$ times.
Among the $n$ tested points,  there is  one  point,  denoted  by
$c_f({\scr P}, n)$,  at which,  $f$ achieves its maximum.
Recall the real maximum of $f$ is achieved  at  $c_f$.

\proclaim{(1.3)Definition} We call $\delta ({\scr P},n):=\sup_{f\in  {\scr F}}\left| c_f-c_f({\scr P},n)\right| $  the
accuracy  of ${\scr P}$ at the $n$-th step.
We say that a policy ${\scr Q}$ is optimal with $n$ steps
if for  any policy ${\scr P}$,  $\delta ({\scr P}, n)\ge  \delta
({\scr Q}, n)$. \endproclaim

\proclaim{(1.4)Theorem (J. Kiefer$^{[3]}$)} The  fraction  method
${{\scr F}}_n$ is optimal with $n$ steps. \endproclaim

At the end of his paper [3],  Kiefer noted that  it  is  not
convenient  in  practice  to use ${{\scr F}}_n$ since  we  have  to
decide in advance the precise number $n$ of tests. Because of this reason,
Kiefer suggested to use $\omega :=\lim_{n\to \infty}F_n/F_{n+1}$
$=(\sqrt  5-1)/2$  as  the first  testing  point  $x_1$
instead  of $F_n/F_{n+1}$ and then keep the  symmetry  rule.  The
later  one  is called the {\it golden section search} (see  also
[8]), denoted by ${\scr W}$. However, L. K. Hua pointed out that ${\scr W}$
is optimal in a different sense and he regarded ${\scr F}_n$ as an approximation
of ${\scr W}$.

\proclaim{(1.5)Definition} A  policy  ${\scr P}$  is   called
symmetric  if  at  the first step,  choose $x_1=x_1({\scr P})$
independent  of  $f\in {\scr F}$.  Starting from the  second  step,
choose  the  new test point according to the  symmetry  rule. \endproclaim

The next result is due to Hua for symmetric policy and extended by J. W. Hong
to the general case.

\proclaim{(1.6)Theorem(Hua$^{[6,7]}$ and Hong$^{[4]}$)}  For any policy  ${\scr
P}$,  we  have $\delta ({\scr P},n)\ge \delta ({\scr W},n) $ for all
sufficiently  large $n$.  In other words,  the policy ${\scr W}$  is
optimal at infinity. \endproclaim

To understand our optimality in a general sense, recall that for each $n$,
the optimal policy  ${{\scr F}}_n$  gives us the  optimal  accuracy
at the last ($n$-th) step:  $\delta ({{\scr F}}_n,n)$$=1/F_{n+1}$.
By  comparing the relative difference (but not the absolute difference
since at different steps, the testing intervals have different scale)
between a policy ${\scr P}$ and the policy ${{\scr F}}_n$
at  the  $n$-th step:
$$ (\delta  ({\scr P},n)-\delta  ({\scr
F}_n,n))/\delta ({{\scr F}}_n,n) =F_{n+1}\delta ({\scr P},n)-1,$$
we arrive at the following notion

\proclaim{(1.7)Definition$^{[1]}$}   We   call
$\delta   ({\scr P}):=\sup_{n\ge  1} F_{n+1} \delta ({\scr P},  n) $
the accuracy of ${\scr P}$.  A policy ${\scr Q}$ is called optimal if for any  policy
${\scr P}$, $\delta ({\scr P})\ge \delta ({\scr Q})$.
\endproclaim

\proclaim{(1.8)Theorem$^{[1]}$} For any symmetric  ${\scr P}$,
we have $ \delta ({\scr P}) \ge \delta ({\scr W})$.
\endproclaim

Before  moving  further,  let us introduce some notations
which will be used throug\-hout this paper.
Let  $k_1,  k_2,  \cdots   (k_1\ge 2) $ be a  sequence  of
positive integers. Set
$$c(k) =  \left(\matrix \format \l &\quad \l \\ \chi (k) &  [\frac
{k+1} 2] \\ \chi (k+1) & [\frac {k+2} 2] \endmatrix \right) $$
where  $[x]$ is the integer part of $x$ and $\chi  (k)=0$  or  $1 $
according to
$k$ being odd or even respectively.  Next, let $ \left(\matrix x_n\\
y_n \endmatrix \right) $ be the solution to the equations
$$ \left(\matrix x_n \\  y_n  \endmatrix \right)  =c(k_{n+1})
\left(\matrix  x_{n+1}  \\  y_{n+1} \endmatrix \right),\quad
n\ge 0,\quad \quad x_0=1.$$
It was proved in [5] that  the
solution  $\left(\matrix x_n \\ y_n \endmatrix \right) $
not only exists but
also unique whenever there are infinitely many odd numbers in the
sequence $\{ k_1, k_2, \cdots  \}. $

We now fix $n\ge 1$ and make the boundary
condition at the final (rather than the first) step:
$ \left(\matrix X_n \\ Y_n
\endmatrix \right) =  \left(\matrix 1 \\ 2 \endmatrix \right)$.
Define  $ \left(\matrix X_m \\  Y_m  \endmatrix \right)
=c(k_{m+1})  \pmatrix  X_{m+1}  \\  Y_{m+1} \endmatrix \right)$,
$0 \le m \le n-1$.
In particular, if $k_m \equiv 2i-1$, then  $X_0=:F^{(i)}_{n+1}$  gives
us  the  generalized Fibonacci sequence$^{[8]}$:
$F^{(i)}_0=F^{(i)}_1=1$ and
$F^{(i)}_n=i(F^{(i)}_{n-1}+F^{(i)}_{n-2})$ for  $n\ge 2$.
For the special case that $k_m \equiv 2i\, (i\ge 1)$, we  rewrite $X_0$ as
$E^{(i)}_n$: $E^{(i)}_n=2(i+1)^n-1$.

\proclaim{(1.9)Definition} Given  an  interval  $[a,b]$  and
$\alpha, \beta    >0$.  The  partition
$a=a_1<b_1=a_2<b_2=a_3<\cdots  <b_N=b $ satisfying $  b_{2k-
1}-a_{2k-1}=\alpha $ and $b_{2k}-a_{2k}=\beta$ for each $k\ge 1$,
if exists, is called an $[\alpha ,\beta]$-partition. \endproclaim

\proclaim{(1.10)Definition}  Fix $n\ge 1$ and $k_1\ge  2,\;  k_2,
\cdots  ,k_n$.  Define the policy ${{\scr G}}_n$ as follows. At
the  first  step,  we take the $[\alpha  ,\beta]$-partition  with
ratio $\alpha /\beta =X_1/(Y_1-X_1)$ and arrange the $k_1$ tests at the
dividing  points.  At  the  $m$-th step,  we  choose  the  $[\alpha
,\beta]$-partition with ratio $\alpha /\beta =X_m/(Y_m-X_m),\; 2\le
m\le n$,  and arrange the $k_m$ new tests, plus the tested point left from
the previous step, at the dividing points. \endproclaim

\proclaim{(1.11)Theorem (Hong$^{[5]}$)}  The policy ${{\scr G}}_n$  is
optimal with  $n$ steps. Moreover, $\delta({{\scr G}}_n,n)$\allowlinebreak
$=1/X_0$.
\endproclaim

\proclaim{(1.12)Definition}  A policy ${\scr P}$ is called  basic
if  at each step the tests are arranged at the dividing points of
an  $[\alpha,   \beta]$-partition.  For  the  special  case  that
$k_n \equiv 2i$  and $\alpha =\beta$ at each step,  the basic policy  is
denoted by ${{\scr E}}^{(i)}$. \endproclaim

The  reason we pay special attention to the basic  policies  is
that  all known optimal policies are basic and on the other hand,
complicated policies are not useful in practice.
Having  these preparations in mind,  we can return to our  main
discussion.  Suppose  that at each step,  $k_m=2i \, (i\ge  1)$.
Then, it is known that
there is no optimal policy at infinity$^{[5]}$.  Nevertheless,  in our new
sense, there does exist an optimal one.

\proclaim{(1.13)Theorem$^{[1]}$} Let $i\ge  2$.  For  any  basic
policy  ${\scr P}$,  we have $\delta ({\scr P}) \ge  \delta  ({\scr
E}^{(i)})=2(i+1)/(2i+1)$. where the accuracy  $\delta
({\scr P})  $  is defined by Definition $(1.7)$ but replacing
$F_{n+1}$ with $E^{(i)}_n $. \endproclaim

The policy ${{\scr E}}^{(i)}$ comes  with  no
surprise since its construction is quite natural.  However,  the case
that  $k_n \equiv 2$ is excluded from the above theorem.  In this  case,
the  optimal policy takes $3/7 $ and $4/7$ as the testing
points at the first step and uses the same construction as  ${\scr
E}^{(i)} (i\ge 2)$ for the subsequent steps $^{[1]}$.

We now consider the odd-block search. That
is,  $k_n \equiv 2i-1\, (i\ge 2).$  In this case,  for fixed $n$, the optimal
policy  ${{\scr G}}_n$ gives us the $[\alpha ,\beta]$-partition:
$\alpha  ={\alpha  }^{(n)}=$$F^{(i)}_n/F^{(i)}_{n+1}$,
$\beta  ={\beta }^{(n)}=1/i-\alpha $ for the first step.  By  the
same   procedure   as  we  mentioned before,
$\omega(i):=$$\lim_{n\to  \infty  }  F^{(i)}_n/F^{(i)}_{n+1}=\big(\sqrt{i(i+4)}-i\big)/(2i)$,
we obtain a basic policy ${{\scr W}}^{(i)}$ by
replacing $\alpha =F^{(i)}_n/F^{(i)}_{n+1} $ with $\omega (i)$ at
the first step.  More precisely,  at the $m$-th $(m\ge 1)$  step, we have
the   $[\alpha   ,\beta]$-partition   with
$\alpha ={\alpha}_m={\omega  (i)}^m$,
$\beta  ={\beta  }_m={\omega  (i)}^{m+1}$.
Moreover,  it  is easy to check that
$\delta ({{\scr W}}^{(i)},n)=
{\omega (i)}^n$,  $n\ge 1$.
Furthermore, it was proved$^{[4,5]}$ in a slight
different  sense that the policy ${{\scr W}}^{(i)}$ is indeed the  optimal
policy at infinity (See Section 4 for details).

Next,  is it true
the policy ${{\scr W}}^{(i)}\, (i\ge 2)$ being the optimal one in  the
general sense? The answer is surprisingly to be negative!

\proclaim{(1.14)Definition} Define  a  basic  policy ${\scr H}={{\scr H}}^{(i)}$
as follows.  At the first step,  we take  the
$[{\alpha }_1, {\beta }_1]$-partition:
$$
{\alpha }_1\!\!=\!\biggl\{
\frac 1 i \biggl[\frac {i+1} 2 \biggr]\!+\!\chi (i)\omega (i)\biggr\}
\omega  (i)
\!=\!\frac  1  i \biggl\{ \chi (i)\!+\!(\chi(i-1)\!-\!\chi(i))\biggl[\frac  {i+1}
2\biggr]\omega  (i) \biggr\}, \quad
{\beta }_1\!=\!\frac 1 i-{\alpha}_1
$$
and at  the  $n$-th  step,  we choose the $[{\alpha  }_n, {\beta  }_n]$-partition:
${\alpha }_n={\omega  (i)}^n$, ${\beta  }_n={\omega (i)}^{n+1}$,
$n\ge 2$. \endproclaim

The  remainder  of  this paper is to prove  the  following result.

\proclaim{(1.15)Theorem} For any basic policy ${\scr P}$, we have
$  \delta  ({\scr P})\ge \delta ({\scr H})$. In  other  words,  the
policy  ${\scr H}={{\scr H}}^{(i)} $ is optimal in the general  sense
among the basic policies. \endproclaim

The  paper is organized as follows.  In the  next  section,  we
prove  some  elementary pro\-perties about the  generalized  Fibonacci
sequence  and a related sequence.  In Section 3,  we study how to
compute the accuracy $\delta ({\scr P},n)$.
At the end of this section,  we explain the main steps of the
proof of Theorem (1.15).  Especially,  we explain why we have  to
study the optimal policy at infinity,  which is the topic studied
in  Section 4.  Having these preparations,  the proof of  Theorem
(1.15) is completed  in  Section 5.  It turns out that the
present proof of the main theorem is quite complicated and lengthy
but we hope that the work would provide some light to solve the problem
for the general situation where the numbers $\{k_n\}$ being arbitrary.

\heading 2. Properties   of   $F$-sequence   and    $G$-sequence
\endheading

  From now on,  we fix $i\ge 2$ and $k_n=2i-1$ at least for all $n\ge
2$.  Thus,  we  may  drop the superscript $i$ from ${{\scr W}}^{(i)},
\; F^{(i)}_n $ and so on without any confusion.

Recall  the  $F$-sequence  is defined  by  $$F_0=F_1=1,\quad\quad F_n=
i(F_{n-1}+F_{n-2}),\quad n\ge 2. \tag 2.1 $$ A related sequence,
called {\it $G$-sequence},  is defined by $$ G_{-1}=0,\quad\quad G_0=1,\quad\quad G_n=i(G_{n-
1}+G_{n-2}),\quad n\ge 1. \tag 2.2 $$

Let $\omega =\omega (i)=\big(\sqrt {i(i+4)}-i\big)/(2i)$ which
is the positive root of $$  i(\omega +{\omega }^2)=1. \tag 2.3 $$

\proclaim{(2.4)Lemma}  For the $F$-sequence, we have
$$\align \text{\hskip-1em}& F_{n+1} =\frac 1 2 \biggl\{ \biggr(1\!+\!3\sqrt {\frac i {i\!+\!4}} \, \biggr)
{\biggl(\frac {i\!+\!\sqrt {i(i\!+\!4)}} 2 \, \biggr)}^n
\!+\!\biggl(1\!-\!3\sqrt {\frac i {i\!+\!4}} \, \biggr){\biggl(\frac {i\!-\!\sqrt
{i(i\!+\!4)}} 2 \, \biggr)}^n \biggr\} \text{\hskip-4em}\tag 2.5\\
&\text{\hskip2.5em}=\frac 1 2 \biggl\{ \biggl(1\!+\!3\sqrt {\frac i
{i\!+\!4}} \, \biggr) {\omega }^{-n}
\!+\!\biggl(1\!-\!3\sqrt {\frac i {i\!+\!4}} \, \biggr){(-i\omega )}^n
\biggr\},\quad \quad n\ge -1.\\
& F_{n+1}F_{n-1}-F^2_n=(2i-1)(-i)^{n-1},\quad\quad n\ge 1. \tag 2.6  \\
& F_nF_{n-1}-F_{n+1}F_{n-2}=(2i-1)(-i)^{n-1},\quad\quad n\ge 2. \tag 2.7 \\
\text{\hskip-1em}&{\text{As}} \; n\to \infty, \; \frac {F_{2n-1}} {F_{2n}} \;
           \text{strictly  increases
          to $\omega $  and } \;
      \frac {F_{2n}} {F_{2n+1}} \; \text{strictly decreases to $\omega$.}
\text{\hskip-3em}\tag2.8 \\
\text{\hskip-1em}&{\text {As}} \; n\to \infty,\; \frac {F_{2n-1}} {F_{2n+1}}\text{ strictly increases
to ${\omega }^2$ and }
  \frac {F_{2n}} {F_{2n+2}} \; \text{strictly decreases to ${\omega }^2$.}
\text{\hskip-3em}\tag2.9 \endalign$$
\endproclaim

\demo{Proof}  Clearly,  (2.5)  follows from (2.1) and  (2.3).
One  may  prove (2.6) by using induction and  (2.1).  Then  (2.7)
follows  from (2.6). Next, by (2.7),  we have
$$\frac {F_{2n+1}} {F_{2n+2}}  -\frac {F_{2n-1}} {F_{2n}}=
\frac {(2i-1)(-i)^{2n}} {F_{2n}F_{2n+2}} >0, \quad\quad
\frac {F_{2n}} {F_{2n+1}}-\frac {F_{2n-2}} {F_{2n-1}}=\frac {(2i-1)(-i)^{2n-1}}
{F_{2n-1}F_{2n+1}} <0.
$$
From this and (2.1),  it is
easy to see that (2.8) holds. Similarly, (2.9) follows from (2.1), (2.7)
and (2.8). \qed \enddemo

\proclaim{(2.10)Lemma}  For the $G$-sequence, we have
$$\align \text{\hskip-1em}& G_n
=\frac 1 {\sqrt {i(i+4)}} \biggl\{ {\biggl(\frac
{i+\sqrt {i(i+4)}} 2 \; \biggr)}^{n+1}-{\biggl(\frac {i-\sqrt {i(i+4)}}
2 \; \biggr)}^{n+1} \biggr\} \tag2.11\\
&\quad \; \, =\frac 1 {\sqrt {i(i+4)}} \{ {\omega }^{-(n+1)}-{(-i\omega )}^{n+1} \},
      \quad\quad n\ge -1. \\
\text{\hskip-1em}&G_nG_m-G_{n+1}G_{m-1}=(-i)^mG_{n-m},\quad\quad n+1\ge m\ge 0.\text{\hskip-3em} \tag 2.12 \\
\text{\hskip-1em}&G_nG_m-G_{n+2}G_{m-2}=-(-i)^mG_{n-m+1},\quad\quad n+1\ge m\ge 1. \text{\hskip-3em}\tag2.13\\
\text{\hskip-1em}&{\text {As}} \; n\to \infty, \; \frac {G_{2n-1}} {G_{2n}}\text{ strictly increases
to $\omega $ and }
\frac {G_{2n}} {G_{2n+1}} \; \text {strictly decreases to} \; \omega.
\text{\hskip-4em}\tag2.14\\
\text{\hskip-1em}&{\text {As}} \; n\to \infty, \; \frac {G_{2n-1}} {G_{2n+1}}\text{ strictly  increases
to ${\omega }^2$ and }\;
\frac {G_{2n}} {G_{2n+2}} \; \text {strictly decreases to} \; {\omega}^2.
\text{\hskip-4em}\tag2.15 \endalign$$
\endproclaim

\demo{Proof}   The  proof  is similar to the  previous  one
except (2.12) and (2.13).  But one may use induction on $m\ge 0$ to prove
that (2.12) holds for all $m \leq n+1 $.  Then,  (2.13) follows from
(2.12) and (2.2). \qed \enddemo

\proclaim{(2.16)Lemma}   Let $a,b,c,d$ be positive numbers  and
$n,m$  be  non-negative  integers with  $n\ge  m$.  Then
$$\frac{aG_m+bG_{m-1}} {aG_{n+1}+bG_n} - \frac {cG_m+dG_{m-1}}
{cG_{n+1}+dG_n}<0 \quad\quad (\text{resp.} >0, \; =0) $$
if and only if $(-1)^m(ad-bc)<0 \; (\text{resp.} >0, \; =0 )$.
\endproclaim

\demo{Proof}  Simply use (2.13) . \qed \enddemo

\proclaim{(2.17)Lemma} For the relation between the $F$-sequence
and $G$-sequence, we have
$$\align
&F_n=G_{n-1}+iG_{n-2},\quad\quad n\ge 1, \tag 2.18 \\
&F_nG_m-F_{n+2}G_{m-2}=-(-i)^mF_{n-m+1},\quad\quad n+1\ge m\ge 1, \tag 2.19\\
&\frac {F_{2n-1}} {F_{2n+1}} < \frac  {G_{2n-
1}} {G_{2n+1}} < \frac {F_{2n+1}} {F_{2n+3}}, \quad\quad
\frac {F_{2n}} {F_{2n+2}} < \frac {G_{2n-2}} {G_{2n}} < \frac {F_{2n-2}}
{F_{2n}}, \quad n\ge 1.
\tag2.20\endalign$$
\endproclaim

\demo{Proof} The first assertion  follows  from  the
definitions of the sequences plus induction.   Then,  (2.19)
follows  from  (2.18) and (2.13).  Finally,  (2.20) follows  from
(2.18) and (2.12). \qed \enddemo

To conclude this section,  we  list the first few terms of the sequences
for the subsequent use.
$$\align &F_0=F_1=1,\quad F_2=2i,\quad F_3=i(2i+1),\quad F_4=i^2(2i+3),
\tag2.21 \\
F_5=i^2(2i^2&+5i+1),\quad F_6=i^3(2i^2+7i+4), \quad
F_7=i^3(2i^3+9i^2+9i+1). \\
&G_{-1}=0,\quad\quad G_0=1,\quad\quad G_1=i,\quad\quad G_2=i(i+1), \tag2.22 \\
G_3=i^2(i+&2),\quad \quad G_4=i^2(i^2+3i+1), \quad\quad G_5=i^3(i^2+4i+3).
\endalign$$

\heading  3.  The Accuracy of ${\scr H}$ and ${\scr W}$.  The Idea of
the Main Proof. \endheading

Suppose that a  policy
${\scr P}$ acts on $f\in {\scr F}$,  after  $n$ steps,  the  remaining
interval  is $[a_n,  b_n]$.  If there is a tested point inside of
$[a_n,  b_n]$,  then let $c_n$  denote this point. Next, set
$$
\Delta ({\scr P},f,n)=b_n-a_n, \quad\quad
\delta ({\scr P},f,n)=
\cases b_n-a_n,\; \text{if there is no $c_n$} \\ \max \{ c_n-a_n,
b_n-c_n \}, \; \text{if $c_n$ exists.}\endcases
$$
For convenience,  put  $\Delta ({\scr P},f,0)=\delta  ({\scr P},f,0)=1$.
Define
$$  \Delta ({\scr P},n)=\sup_f \Delta ({\scr P},f,n),\qquad
\delta ({\scr P},n)=\sup_f \delta ({\scr P},f,n).$$
The last one is the accuracy of ${\scr P}$ at the $n$-th step,
which is precisely the same as we defined before.  Finally,  the
accuracy  of ${\scr P}$ is given by
$ \delta ({\scr P})=\sup_{n\ge 1}  F_{n+1}  \delta ({\scr P},n)$,
since we are now in  the  case that $k_n=2i-1$.

As we have seen in the first section, $\delta({\scr W}, n)={\omega }^n,\;
n\ge 0 $. We now prove

\proclaim{(3.1)Lemma} $\delta ({\scr W})=F_2 \omega =2 i \omega .$
\endproclaim

\demo{Proof} The proof will be done once we show that  as
$n\to  \infty, \; F_{2n+1}{\omega  }^{2n}  $   increases   and
$F_{2n+2}{\omega }^{2n+1}$ decreases to the same limit
$$\lim_{n\to \infty } F_{n+1}{\omega }^n =\frac 1 2
+\frac 3 2 \sqrt {\frac i {i+4}} =\frac 1 {i+4}
\big(2(i+1)+3i\omega \big). \tag3.2$$
The last conclusion follows from (2.5) immediately. On the other
hand, by (2.9), we have
$$\frac {F_{2n+3}{\omega }^{2n+2}} {F_{2n+1}{\omega }^{2n}}
= {\omega }^2\bigg /\frac {F_{2n+1}} {F_{2n+3}} \;  \downarrow  1,
\quad\quad \frac {F_{2n+2}{\omega }^{2n+1}} {F_{2n}{\omega }^{2n-1}} =
{\omega }^2\bigg/ \frac {F_{2n}} {F_{2n+2}} \; \uparrow  1.
$$
Hence the proof is completed. \qed \enddemo

As for the policy ${\scr H}$ defined in (1.14), we have

\proclaim{(3.3)Lemma}  $\delta ({\scr H})=F_4{\omega}^3=i^2(2i+3){\omega }^3
=: \delta $. \endproclaim

\demo{Proof} The last equality follows from (2.21).  By (2.3), we can
also express $\delta$ as follows:
$$ \align  \delta  &=i(2i+3)(1-i\omega   )\omega
=i(2i+3)((i+1)\omega  -1) \tag 3.4 \\ &=\frac {2i+3} 2 \bigl((i+1)\sqrt
{i(i+4)}-i(i+3)\bigr) <2. \endalign  $$
By the  definition of ${\scr H}$,  it is easy to check that
$$  \delta
({\scr H},1)=\biggl\{ \frac 1 i \biggl [\frac {i+1} 2\biggr]+
\chi (i)\omega \biggr \} \omega, \quad \quad
\Delta ({\scr H},1)=1/i $$ and
$\delta ({\scr H},n)={\omega  }^n$,
$\Delta ({\scr H},n)={\omega }^{n-1}/i$,  $n\ge 2$.
Thus,  it follows from  the  proof of Lemma (3.1) that
$$ \align
\delta ({\scr H})&= \sup_{n\ge
1}  F_{n+1}\delta ({\scr H},n) \\
&=\max \Big\{F_2\delta  ({\scr H},1),\; \sup_{n\ge  2}F_{n+1}\delta ({\scr H},n)\Big\} \\
&= \max \{ F_2\delta ({\scr H},1),\; F_4{\omega }^3 \}.
\endalign
$$
Hence, we need only to show that
$$2i\biggl (\frac 1 i\biggl [\frac {i+1} 2\biggr ]+\chi (i)\omega
\biggr )<i(2i+3)(1-i\omega).\tag 3.5 $$
If $i=2$,  then $\omega =(\sqrt 3-1)/2$. In this case, a direct
computation shows that (3.5) holds.  For $i\ge 3$,  the left hand
side  of  (3.5) is less or equal to
$i+1+2i\omega$.
Thus,  it  suffices  to  show  that
$i(2i^2+3i+2)\omega   <2i^2+2i-1$. Noticing   that   $\omega
<F_4/F_5=(2i+3)/(2i^2+5i+1),$   the  above   inequality   follows
immediately for all $i\ge 3$. \qed \enddemo

\proclaim{(3.6)Corollary} $\delta ({\scr W}) > \delta ({\scr H}).$
\endproclaim

This corollary means that ${\scr H}$ is better than ${\scr W}$ in  the
general sense. Comparing ${\scr H}$ with ${\scr W}$ carefully, we see
that  the difference between these two policies is only  at
the  first step.  For ${\scr H}$,  we choose
$${\alpha }_1({\scr H})
=\biggl \{ \frac 1 i\biggl [\frac {i+1} 2\biggr ]+\chi (i)\omega
\biggr \} \omega, \quad\quad {\beta }_1({\scr H})=1/i-{\alpha }_1({\scr H}).
$$
For ${\scr W}$, we choose
${\alpha }_1({\scr W})=\omega $, ${\beta }_1({\scr W})=1/i-{\alpha }_1({\scr W})$.
Starting
from  the second step,  the construction rule for the two policies
is completely the same. Now, what is the key point to making such
a  choice for ${\scr H}$.  The reason is as follows.  Since at  the
first step, we have an odd-block search $k_1=2i-1$, for any basic
policy ${\scr P}$,  we always have
$\Delta ({\scr P},1)=1/i$. On
the  other hand,  at the second step,  including the tested point
(left  from the first step),  there are altogether  $2i$  tests.
However, the $\ell$-th position $(1\le \ell \le 2i)$ located by the tested
point  left  from the first step does not make any influences  to
the  construction  for  the second step,  since the  key  of  the
construction is the ratio ${\alpha }_2/{\beta }_2 $.  This is due
to the fact that our policies are basic.  But the location of the
tested  point does influence $\delta  ({\scr P},1)$.  Furthermore,
each  basic  policy ${\scr P}$ corresponds uniquely a  basic  policy
with  initial  testing  interval  $[0,  1/i]$  and  with  testing
numbers:  $k_1=2i,\;  k_n=2i-1,\;  n\ge  2$.  Let us denote the later
policy by ${{\scr P}}_1$. Corresponding to ${\scr H}$, we have ${\scr
H}_1$.  Conversely,  due  to the rule for the basic  policies,  a
basic  policy ${{\scr P}}_1$ with initial testing
interval  $[0,1/i]$ and with testing numbers  $k_1=2i,\;  k_n=2i-1,\;
n\ge 2$, determines uniquely (here we regard those policies which have
the same accuracy at the first step as the same) a basic policy
${\scr P}$ with initial
testing  interval  $[0,1]$ and with testing  numbers  $k_n\equiv 2i-1$.
Moreover, It is obvious that
$$ \delta ({\scr P})=\sup_{n\ge 1} F_{n+1}\delta ({\scr P},n)
=\sup_{n\ge 1} F_{n+2}\delta ({{\scr P}}_1,n). \tag 3.7 $$
In particular,
$$\delta ({\scr H})=\sup_{n\ge 1} F_{n+2}\delta ({{\scr H}}_1,n)
=F_4{\omega }^3. \tag 3.8 $$
The above discussions tell  us,  in
order  to prove Theorem (1.15),  we need only to study the  basic
polices ${{\scr P}}_1$ and proving that
$$\sup_{n\ge 1}F_{n+2}\delta ({{\scr P}}_1,n)\ge  F_4{\omega
}^3=\sup_{n\ge 1} F_{n+2}\delta ({{\scr H}}_1,n). \tag 3.9 $$
To  fix  our idea,  let us repeat the construction of  the  basic
policy  ${{\scr H}}_1$.  At the $n$-th step,  we take the  $[{\alpha
}_n,{\beta  }_n]$-partition with
${\alpha  }_n={\omega  }^{n+1}$,
${\beta }_n={\omega }^n$.

\proclaim{(3.10)Definition}  We say that two policies ${\scr P}$
and  ${{\scr P}}'$ are equivalent if for all  $n\ge 1$,
$\delta ({\scr P},n)= \delta ({{\scr P}}',n)$.
\endproclaim

Again,  due to the rule  of
the basic policies, it is easy to check the following fact:

\proclaim{(3.11)Lemma}  If  a basic policy ${{\scr P}}_1$ is not
equivalent  to  ${{\scr H}}_1$,  then  we  must  have
$  {\alpha
}_1={\alpha }_1({{\scr P}}_1)\ne {\omega }^2= {\alpha }_1({\scr
H}_1)$.
\endproclaim

We now going to study how to compute the accuracy at each step.

\proclaim{(3.12)Remark}  Let $[a_n,c_n,b_n]$ be  the  remaining
interval  left from the $n$-th step,  for the purpose of  computing
the accuracy at the $(n+1)$-th step, we may and will assume that
$a_n=0$, $c_n=\delta ({{\scr P}}_1,n)$, $b_n=\Delta ({{\scr P}}_1,n)$
and $c_n\ge b_n/2 $. Moreover, since $k_n=2i-1$ for all $n\ge 2$,
we may also assume that $c_n>b_n/2$. \endproclaim

\proclaim{(3.13)Definition}  Let $[0,\; \delta ({{\scr P}}_1,n),
\; \Delta ({{\scr P}}_1,n)]$ be the remaining interval from the n-th
step, $n\ge 0$. For simplicity,  we write  ${\delta }_n=\delta
({{\scr P}}_1,n)$ and ${\Delta }_n=\Delta ({{\scr P}}_1,n)$
respectively.   At  the  $(n+1)$-th step, there are $2i$ tests
including  the  tested point ${\delta }_n$.  Denoted them  by
$0=z_0<z_1<\cdots <z_{2i}<z_{2i+1}={\Delta }_n$.
Then,  it
follows  from  the last remark that ${\delta }_n$ must be one  of
$\{ z_{i+1}, z_{i+2},\cdots  ,z_{2i}  \}$.  If  ${\delta
}_n=z_{\ell}$,  we  say  that  ${\delta }_n$ is located  at  the $\ell$-th
position. \endproclaim

\proclaim{(3.14)Lemma} Let ${\delta }_{m-1}$ be located at the
$\ell$-th position. Then
$$\left(\matrix {\alpha }_m \\ {\Delta }_m \endmatrix \right)
=\frac 1 {i\,\chi (\ell-1)-[\ell/2]} \left(\matrix i{\delta
}_{m-1}-[\ell /2]{\Delta}_{m-1} \\ \chi (\ell-1){\Delta }_{m-1}
-{\delta  }_{m-1} \endmatrix \right) \tag 3.15 $$
This occurs only if
$${\delta }_{m-1}>\frac {\ell } {2(i+1)} {\Delta }_{m-
1},\; \text{if $\ell=2j$,} \quad\quad {\delta }_{m-1}
< \frac {\ell +1} {2(i+1)}
 {\Delta }_{m-1},\; \text{if $\ell =2j-1$.}
\tag3.16 $$
Then, ${\delta }_m={\alpha }_m$  if and only if
$${\delta }_{m-1}<\frac {\ell } {2i+1} {\Delta }_{m-1},\;
\text{if $\ell =2j$,} \quad\quad {\delta }_{m-1}>\frac {\ell } {2i+1}
{\Delta }_{m-1},\; \text{if $\ell =2j-1$.}
\tag3.17 $$
Otherwise,
$${\delta }_m = {\Delta }_m-{\alpha }_m =
\frac {(\chi(\ell -1)+[\ell /2]){\Delta }_{m-1}-(i+1){\delta
}_{m-1}} {i\chi (\ell -1)-[\ell /2]}. \tag 3.18
$$ \endproclaim

\demo{Proof}  Since ${\delta }_{m-1} $ located at the $\ell $-th
position, we have
$$\left(\matrix {\delta }_{m-1} \\ {\Delta }_{m-1} \endmatrix \right)
=\left(\matrix \format \l  &\quad \l \\ \chi (\ell -1)
&[\ell /2] \\ 1 & i \endmatrix \right)  \left(\matrix  {\alpha  }_m  \\
{\Delta }_m \endmatrix \right) \tag 3.19 $$ and
$$ {\Delta }_m > {\alpha }_m. \tag 3.20 $$
Now,  (3.15) follows from (3.19).  By using (3.15), it is easy to
check  that  (3.16) is equivalent to (3.20).  On the other  hand,
${\delta  }_m={\alpha }_m$ is equivalent to say that
$  {\alpha}_m >{\Delta }_m/2 $ by Remark (3.12).  Hence,  the last
two assertions follows by a simple computation. \qed \enddemo

Let
$ j_e= i/2 +1 $ and $j_0=(i+1)/ 2 $
for even  and odd $i$ respectively.  Define
$ K(j)=  j {\Delta  }_{m-1}/(i+1)$.

\proclaim{(3.21)Lemma} There are altogether four cases:

$i)$ If ${\delta }_{m-1} \in (K(j), K(j+1)) $ and ${\delta }_{m-1} $
is located at an even's position $\ell $,  then $\ell =2j$.  Where j varies
from $j_e$ to $i-1$ if $i$ is even; otherwise, from $j_0$ to $i-1$.

$ii)$  If ${\delta }_{m-1}\in (K(j),K(j+1))$ and ${\delta }_{m-1}$  is
located  at an odd's position $\ell $,  then $\ell =2j+1$.  Where $j$  varies
from $j_e$ to $i-1$ if $i$ is even; otherwise, from $j_0+1$ to $i-1$.

$iii)$ If ${\delta }_{m-1}\in (K(i),  {\Delta }_{m-1})$, then ${\delta
}_{m-1}$ can only be located at the $(2i)$-th position.

$iv)$  If ${\delta }_{m-1}\in ( {\Delta  }_{m-1}/2,  K([(i+3)/2]))$,
then ${\delta }_{m-1}$ can only be located at  the
${\ell }_0$-th position with
${\ell }_0=2\bigl [\frac {i+3} 2\bigr ]-1=i+1-\chi (i+1)$. \endproclaim

\demo{Proof} The partition of sub-intervals $(K(j),K(j+1))$
is suggested by Lemma (3.14), especially (3.16). The range of $j$ 's
is  due to the fact that the position $\ell $ varies from $i+1$ to $2i$ and
the  fact  that $\frac 1 2 {\Delta }_{m-1} <  {\delta  }_{m-1}  <
{\Delta }_{m-1}.$ \qed \enddemo

\proclaim{(3.22)Lemma}   Given  a  basic policy  ${{\scr P}}_1$,
assume  that at the $m$-th step,  we have the $[{\alpha  }_m,{\beta
}_m]$-partition.

$i)$ If ${\alpha }_1\in ({F_{2n-1}}/{F_{2n+1}}, {G_{2n-2}}/{G_{2n}})$, then
$$\delta ({{\scr P}}_1,m)=(-1)^m (G_{m-2}-G_m {\alpha }_1)/i^m,
\quad\quad 1\le m\le 2n-1 \tag 3.23 $$
$$\Delta ({{\scr P}}_1,m)=\delta ({{\scr P}}_1,m-1)/i, \quad\quad 2\le
m\le 2n-1. \tag 3.24 $$
Moreover,  $\delta  ({{\scr P}}_1,m-1)={\alpha }_m >{\beta }_m$ for
all $1\le m\le 2n-1$.

$ii)$  If  ${\alpha  }_1\in  ({G_{2n-1}}/{G_{2n+1}}, {F_{2n}}/{F_{2n+2}})$,
then the same assertions in $i)$ hold  for
all $m$ up to $2n$.  \endproclaim

\proclaim{(3.25)Remark} If we set $G_{-2}=1/i$,  then we  can
keep not only the recurrence $G_m=i(G_{m-1}+G_{m-2}),\; m\ge 0$ but
also extend $(3.24)$ to $m=1$,  regarding $\delta ({{\scr P}}_1,0)  $
as those given by the right hand side of $(3.23)$.  We will use this
convention for simplicity.  However,  we will use this convention
only  for  computing  ${\delta }_m$'s and ${\Delta }_m$'s  with
starting  value  $m=1$.  Otherwise,  it would contradict  to  our
original  convention  that $\delta  ({{\scr P}}_1,0)=1/i$  which  is  the
length of the initial testing interval. \endproclaim

\demo{Proof of Lemma $(3.22)$}

a) Let us begin with the first step. We have ${\Delta }_1=\frac 1 i
\bigl (\frac 1 i-{\alpha }_1 \bigr ) > {\alpha }_1$. This gives  us  one
condition
${\alpha }_1 < \frac 1 {i(i+1)} = \frac {G_0} {G_2}$.
Clearly,  ${\delta }_1={\alpha }_1$ if and only if
$ {\alpha }_1 >  \frac  1 {i(2i+1)} =\frac {F_1} {F_3}$.
Hence ${\alpha  }_1\in (\frac {F_1} {F_3},\frac  {G_0} {G_2} )$.

b) Next, consider the second step. Assume that ${\alpha }_1\in
(\frac {G_1} {G_3}, \frac {F_2} {F_4})\subset (\frac {F_1} {F_3},
\frac {G_0} {G_2})$.   Since  in  the  present  situation,   $
K(i)=\frac i {i+1}{\Delta }_1$,  by a),  it  follows  that
${\delta }_1\in (K(i),  {\Delta }_1)\Leftrightarrow {\alpha }_1\in
\bigl(\frac  {G_1}  {G_3},  \frac {G_0} {G_2}\bigr) $.
But
$ \bigl(\frac {G_1} {G_3},
\frac {F_2} {F_4}\bigr) \subset \bigl(\frac {G_1} {G_3},
\frac {G_0} {G_2}\bigr)$,
hence by Lemma (3.21) we see that ${\delta }_1$ must be
located   at  the  $(2i)$-th  position.   On the  other   hand,   by
(3.17), ${\delta }_2={\alpha }_2$ is equivalent to ${\alpha }_1 <
F_2/F_4$. Thus, we have proved the lemma for $n=1$.

c) Suppose that we now arrive at the $m$-th step.  Then  ${\delta
}_m > K(i) $ means that
$${\alpha }_m >  i {\Delta }_m/(i+1) = {\delta }_{m-1}/(i+1). \tag 3.26 $$
Noticing that we have already had ${\delta }_m < {\Delta }_m$  by
our  assumption,  hence  by Lemma (3.21),  ${\delta }_m$ must  be
located  at the $(2i)$-th position whenever (3.26)  holds.  On  the
other hand,  by (3.17),  ${\delta }_{m+1}={\alpha }_{m+1}$ if and
only if
$$(2i+1){\delta }_m < 2{\delta }_{m-1}. \tag 3.27 $$
Furthermore,  if these all hold, then
$ {\delta }_{m+1}= {\alpha}_{m+1}  ={\delta }_{m-1}/i-{\delta }_m$.
By the hypotheses of
induction,  we then obtain
$$\align
 {\delta }_{m+1} &=\frac  {(-
1)^{m-1}} {i^{m-1}} (G_{m-3}-G_{m-1}{\alpha }_1) - \frac {(-1)^m}
{i^m}(G_{m-2}-G_m{\alpha  }_1)\\
& = \frac {(-1)^{m+1}} {i^{m+1}}
(G_{m-1}-G_{m+1}{\alpha }_1).
\endalign$$
This gives (3.23).

We now return to (3.26).  By the hypotheses of induction, (3.26)
is equivalent to
$$(i+1)(-1)^m  (G_{m-2}-G_m{\alpha  }_1)/i^m >
(-1)^{m-1} (G_{m-3}-G_{m-1}{\alpha }_1)/i^{m-1}. \tag 3.28 $$
For odd $m$,  this becomes
$$ G_{m-3}+(i+1) G_{m-2}/i < (G_{m-1}+ (i+1) G_m/i){\alpha }_1.$$
By  (2.2),  this  gives  us
${\alpha }_1>G_m/G_{m+2}$. For even $m$,  (3.28) is equivalent to
${\alpha }_1 < G_m/G_{m+2}$. On the other hand,  for odd  and
even $m$, (3.27) gives us
$$ {\alpha }_1 < (2G_{m-1}+G_{m-2})/(2G_{m+1}+G_m)
=F_{m+1}/F_{m+3}$$
and
${\alpha }_1  >
F_{m+1}/F_{m+3}$  respectively. Combining these facts, we prove
the required conclusions. \qed \enddemo

Observing (2.9), (2.15),  (2.20), Lemma (3.11) and Lemma (3.22),
it is natural to assume that ${\alpha }_1={\alpha }_1({{\scr P}}_1)
$ is in one of the following sub-intervals:
$$\align \text{\hskip-2em}&\biggl (0,\frac {F_1} {F_3}\biggr ),\quad\quad \biggl (\frac {F_{2n-1}}
{F_{2n+1}},\frac {G_{2n-1}}  {G_{2n+1}} \biggr ), \quad\quad \biggl (\frac
{G_{2n-1}} {G_{2n+1}}, \frac {F_{2n+1}} {F_{2n+3}}\biggr )\text{\hskip-5em}\tag3.29 \\
&\biggl (\frac {G_{2n}} {G_{2n+2}}, \frac {F_{2n}} {F_{2n+2}} \biggr ),
\quad\quad \biggl (\frac {F_{2n}} {F_{2n+2}},  \frac {G_{2n-2}}
{G_{2n}}\biggr ), \quad\quad \biggl (\frac {G_0} {G_2}, \frac 1 i \biggr ),
\quad \quad n\ge 1. \endalign $$
The first three are contained in $[0,{\omega }^2]$ but the second
three in $[{\omega }^2, 1/i]$. However, as we have seen from the
proof a) of Lemma (3.22), $ {\alpha }_1 < G_0/G_2$. Thus,  the
last one in (3.29) can be ignored. Now, we want to prove that the
second  one for $n\ge 2$ and  the third one for all $n$ can also be
ignored.

\proclaim{(3.30)Lemma}  Let  $m\ge 1$ and
$ {\delta}_m =(-1)^m (G_{m-2}-G_m{\alpha}_1)/i^m$.
Then  $F_{m+2}{\delta}_m < \delta $ if and only if
$$(-1)^{m-1} {\alpha }_1 < \bigl\{ i^m\delta/F_{m+2} +
(-1)^{m-1} G_{m-2} \bigr\}/G_m=: (-1)^{m-1} A_m. \tag 3.31 $$
In particular,
$ A_1=\delta /F_3=F_4{\omega }^3/F_3$ and $A_2={\omega}^2$.
\endproclaim

\demo{Proof}    The  assertions  follow  by   some   simple
computations.  For  instance,  let $m=2$,  then
$${\alpha }_1 > \frac 1 G_2 \biggl( G_0-\frac {i^2\delta } {F_4} \biggr)
=\frac 1 {i(i+1)}  \{  1-i^2{\omega  }^3 \} =\frac 1  {i(i+1)}  \{  1-
i\omega (1-i\omega ) \} ={\omega }^2.
$$
This shows that $A_2={\omega }^2$. \qed \enddemo

By   Lemma  (3.22), if ${\alpha }_1\in (\frac {F_{2n-1}}
{F_{2n+1}}, \frac {G_{2n-1}} {G_{2n+1}})$ with $n\ge 2$ or ${\alpha
}_1\in (\frac {G_{2n-1}} {G_{2n+1}}, \frac {F_{2n+1}} {F_{2n+3}})$
with $n\ge 1$,  then the assumption of Lemma (3.30) holds for $m=1$
and $2$.  Hence $F_4{\delta }_2 \ge \delta$ and so (3.9)  holds.
Thus,  the  proof  of our main theorem is done for  these  cases.
Therefore,  we need only to consider the cases that ${\alpha }_1$
is in one of the sub-intervals:
$$ \biggl (0,\frac {F_1} {F_3}\biggr ),\quad\quad \biggl (\frac {F_1} {F_3},
\frac {G_1} {G_3} \biggr ),\quad\quad \biggl (\frac {G_{2n}} {G_{2n+2}},
\frac {F_{2n}} {F_{2n+2}} \biggr),\quad\quad \biggl (\frac {F_{2n}} {F_{2n+2}},
\frac {G_{2n-2}} {G_{2n}} \biggr ). \tag 3.32 $$
Suppose  that ${\alpha }_1$ is in the second interval.  Then,  by
Lemma (3.22), the explicit expression for $\delta ({\scr
P}_1,m)$ works only for $m=1$. If we could find an $n_0$ so that
$$ F_{n_0+2}\delta ({{\scr P}}_1, n_0) \ge \delta. \tag 3.33 $$
Then,  we were done. The problem is that such an $n_0$ for which (3.33) holds
can be very large. And there is no simple way to find out $n_0$ since
there is no simple expression for $\delta  ({\scr
P}_1,m)$ when $m>1$. Because
of  this,  we employ the limiting behavior of $F_{n+2}\delta
({{\scr P}}_1,n)$ as $n\to \infty$.  And this is just what  we  are
going to study in the next section.

\heading 4. Optimal Policy at Infinity. \endheading

In this section,  we study the optimal policy at infinity.  The
results  obtained  here are not only for the later use  but  also
have their own interesting.
For  our  reader's convenience,  we first copy some lemmas
from  [5]  which are available for any  sequence  $k_1\ge  2,\; k_2,\;
\cdots $ of positive integers.
Let $ \left(\matrix x_n \\ y_n \endmatrix \right) $ be the
solution to the equations
$$\left(\matrix x_n  \\ y_n \endmatrix \right)
= c(k_{n+1}) \left(\matrix x_{n+1} \\ y_{n+1} \endmatrix \right),
\;\; n\ge 0,\quad\quad  x_0 =b-a>0, \tag4.1 $$
and  $\left(\matrix u_n \\ v_n  \endmatrix \right) $ satisfies
$$\left(\matrix  u_n  \\  v_n   \endmatrix \right)  \le
c(k_{n+1}) \left(\matrix u_{n+1} \\ v_{n+1} \endmatrix \right),
\;\; n\ge 0, \quad\quad u_0=b-a. \tag4.2 $$
Put
$$ \mu (m,n)=\frac {v_n} {y_n}\bigg/\frac {u_m} {x_m},\quad\quad
\lambda (m,n)= \frac {v_n} {y_n}\bigg/\frac {v_m} {y_m},\quad\quad
\rho (m,n) =\frac {u_n} {x_n}\bigg/\frac {u_m} {x_m}.  \tag4.3$$
Obviously, we have

\proclaim{(4.4)Lemma$^{[5]}$} \;\;  $\lambda (m,l)\lambda (l,n)=\lambda (m,n)$,
\quad $\rho (m,l)\rho (l,n)=\rho (m,n)$,  \newline
\text{\hskip8.0em} $\mu (m,l)\lambda  (l,n)=\mu (m,n)$, \quad
$\rho (m,l)\mu (l,n)=\mu (m,n)$, $m,n,l \ge 0$. \endproclaim

\proclaim{(4.5)Lemma$^{[5]}$} If $k_{n+1}=2i$, then $\lambda (n,n+1)
\ge 1$ and $\mu (n,n+1)\ge i/(i+1)$.
If  $k_{n+1}=2i-1$,
then $ \lambda (n,n+1)\ge i/(i+1)$ and  $\mu (n,n+1)\ge 1$.
\endproclaim

\proclaim{(4.6)Lemma$^{[5]}$} If $\lambda (n,n+1)<1$, then $\mu (n+1,n)^{-1}
\ge \lambda (n,n+1)^{-1}$. \endproclaim

\proclaim{(4.7)Lemma$^{[5]}$} If $ \mu (n,n+1)<1$,  then $\rho  (n,n+1)
\ge \mu (n,n+1)^{-1}$. \endproclaim

Now, we return to our main setup: $k_1\ge 2,\; k_n=2i-1 $ for all
$ n\ge 2$. In this case, by Lemma (4.5), we have
$$\lambda (n,n+1)\ge i\dsize /(i+1) \; \text{and $\mu (n,n+1)\ge 1$ for all
$n\ge 2$}. \tag 4.8 $$

\proclaim{(4.9)Lemma} For each $n\ge 0$,  if $ \lambda  (n,n+1)
<1$, then $ \lambda (n,n+2)\ge \lambda (n,n+1)^{-1}  > 1$.
\endproclaim

\demo{Proof}  By  Lemma  (4.4),  it  follows  that  $\mu
(n+1,n)\lambda (n,n+2)=\mu (n+1,n+2)$. Hence, by Lemma (4.6), we
have
$\lambda (n,n+2)=\mu(n+1,n)^{-1}\mu (n+1,n+2)
\ge \mu (n+1,n)^{-1} \ge \lambda (n,n+1)^{-1}$.\qed\enddemo

\proclaim{(4.10)Lemma$^{[5]}$} Let
$   \sigma  =(b-a) \big /\bigl \{ \chi (k_1) \omega +\bigl [\frac {k_1+1} 2\bigr ]\big /i \bigr \}$.
Then
$x_n=\sigma{\omega}^n$, $y_n=\sigma {\omega }^{n-1}/i$, $n\ge 1$  is  the
unique   solution  to  the  equations
$$  \left(\matrix x_n \\ y_n \endmatrix  \right)
=  c(k_{n+1})  \left(\matrix  x_{n+1}\\y_{n+1} \endmatrix \right),
\; n\ge 0,\quad\quad x_0=b-a. $$
\endproclaim

From now on,  unless otherwise stated,  let ${{\scr P}}_2$  to
denote  an arbitrary policy with initial testing interval $[a,b]$
and with successive testing numbers $k_1 \ge 2, \; k_n=2i-1$ for all
$n\ge 2$.  Let ${{\scr H}}_2$ denote the basic policy:  at the  $n$-th
step,  the $[{\alpha }_n, {\beta }_n]$-partition is determined by
$ {\alpha  }_n=x_n,\; \, {\beta  }_n=y_n-x_n,\;\, n\ge 1$,
where
$(x_n,y_n)$ is given by Lemma (4.10). For simplicity, we put
$ u_n=\delta({{\scr P}}_2,n)$, $v_n=\Delta ({{\scr P}}_2,n)$, $n\ge 0 $ where
$u_0=\delta
({{\scr P}}_2,0)=b-a$. It is known that $\left(\matrix u_n \\ v_n
\endmatrix \right) $ satisfies $(4.2)^{[5]}$.

Actually, the policies ${{\scr P}}_2$ and ${{\scr H}}_2$  are  the
generalization of ${{\scr P}}_1$ and ${{\scr H}}_1$ respectively. If we take $[a,b]=[0,1/i]$
and $k_1=2i$, then ${{\scr P}}_2$ and ${{\scr H}}_2$
coincide with  ${{\scr P}}_1$ and ${{\scr H}}_1$ respectively.  But
we prefer to distinguish them.

Next, we introduce  a sequence $\{ {\varphi }_j \}  $  by  the
following procedure.

(I) Let $k_1$ be an odd integer.  In this case,  we always have
$\mu (0,1)\ge 1$. If $\mu (0,1)>1$, we simply take ${\varphi }_1=\mu (0,1)$.
Otherwise,  we  look at the sequence  $\{  \lambda
(n,n+1):  n\ge 1 \}$. If $ \lambda (n,n+1)\ge 1 $ for all $n\ge
1$,  we  cancel  those $\lambda (n, n+1) $ which equals one  and
denote  by  $  {\varphi }_1,  {\varphi }_2,  \cdots  $  the
remaining   $\lambda   (n,n+1)$'s successively. Then the
construction  is  done.  Conversely,  if we find some  $  \lambda
(n,n+1)<1$.  Then,  by Lemma (4.9), we have $ \lambda (n,n+2)\ge
\lambda (n,n+1)^{-1} > 1$. In this case, we will forget $\lambda
(n,n+1)$ and $ \lambda (n+1,n+2)$ and take $\lambda (n,n+2) $  as
one  of the $\varphi$'s.  Then go ahead to look at $\{  \lambda
(m,m+1):  m\ge n+2 \} $ and repeat the same procedure. Of course,
the  index  set  $J$ of $ \{ {\varphi }_j \} $ may be  empty,  which  is
equivalently to say that $ \mu (0,1)=1 $ and $ \lambda  (n,n+1)=1
$ for all $n\ge 1$.

(II) Let $k_1$ be an even integer. If $\mu (0,1)\ge 1$, then we
can  adopt the same construction for $\varphi $ as given in  (I).
But  in this case,  it can be happen that $\mu (0,1)<1$ for which
we have to modify the above construction. By Lemma (4.7), we then
have $ \rho (0,1)\ge \mu(0,1)^{-1} >1$.
And so we set ${\varphi
}_1=\rho (0,1)>1$. Now, we have $\mu (1,2)\ge 1$. This enables us
to  return to the previous construction by regarding $\mu  (1,2)$
and  $\{  \lambda (n,n+1):  n\ge 2 \} $ as $\mu  (0,1)$  and  $\{
\lambda (n,n+1):  n\ge 1 \} $ respectively.  Again, the index set
$J$ of  $\{  {\varphi  }_j  \}  $  is empty if  and  only  if  $  \mu
(0,1)=\lambda (n,n+1)=1 $ for all $n\ge 1$.

To get  a precise impression of the  above  construction,
consider a special case:
$$ \mu (0,1)<1, \; \mu(1,2)>1, \;\;
         \lambda (2,3)=1, \; \lambda (3,4)<1
\;\text{ and }\; \lambda (n,n+1)=1 \; \text{for all} \;n\ge 5.
$$
Then we have
$ {\varphi  }_1=\rho (0,1)$, ${\varphi  }_2=\mu(1,2)$,
${\varphi}_3=\lambda (3,5)$
and $J=\{ 1,2,3 \}$. On the other hand,  by Lemma
(4.4), we see that
$$\align
\frac {\Delta ({{\scr P}}_2,n)} {\Delta
({{\scr H}}_2,n)} &=\frac {v_n} {y_n}=\frac {v_n} {y_n}\bigg /\frac {u_0}
{x_0} =\mu (0,n) =\rho (0,1)\mu (1,n)
=\rho (0,1)\mu (1,2)\lambda (2,n)   \\
&=[\rho (0,1)] \; [\mu (1,2)] \; \lambda (2,3) \;
  [\lambda (3,4)\lambda (4,5)] \; \lambda (5,6)\cdots  \lambda (n-1,n)
\ge \dsize \prod_{j=1}^3 {\varphi }_j,\\
&\qquad\qquad n\ge 5.
\endalign  $$
This  example not only shows the  reason  why  we
introduced  such  a  construction for ${\varphi }'$s  but  also
indicates the proof of the following result.

\proclaim{(4.11)Proposition$^{[4,5]}$} For any ${{\scr P}}_2$,  we
have
$\dsize\lim_{n\to  \infty }  \dfrac  {\Delta  ({{\scr P}}_2,n)}
{\Delta  ({{\scr H}}_2,n)} \ge \dsize \prod_j {\varphi }_j \ge 1$.
Moreover,  the equality holds if and only if $\mu  (0,1)=\lambda
(n,n+1)=1$ for all $n\ge 1$. \endproclaim

We are now ready to prove the main result of this section.

\proclaim{(4.12)Theorem}  ${{\scr H}}_2$ is the optimal policy at
infinity.  That  is,  for  any policy ${{\scr P}}_2$, we have $$
\lim_{n\to \infty } \frac {\delta ({{\scr P}}_2,n)} {\delta  ({\scr
H}_2,n)}  =\lim_{n\to  \infty  } \frac  {\Delta  ({{\scr P}}_2,n)}
{\Delta ({{\scr H}}_2,n)} = \dsize \prod_j {\varphi }_j \ge 1. $$
\endproclaim

\demo{Proof}  Write  $\delta ({{\scr P}}_2,n)/\delta  ({\scr
H}_2,n) $ as
$$ \frac {\delta ({{\scr P}}_2,n)} {\delta ({{\scr H}}_2,n)}
=\frac {\Delta ({{\scr P}}_2,n)} {\Delta ({{\scr H}}_2,n)}\cdot
\frac {\delta ({{\scr P}}_2,n)} {\Delta ({{\scr P}}_2,n)i\omega }
=D_n\frac {{\xi }_n} {i\omega} \; ,
\tag 4.13 $$
where   $$D_n=\frac {\Delta ({{\scr P}}_2,n)} {\Delta
({{\scr H}}_2,n)},\quad\quad {\xi }_n  = \frac {\delta ({{\scr P}}_2,n)}
{\Delta ({{\scr P}}_2,n)} . $$
Here, we have used the fact that $
\delta ({{\scr H}}_2,n)=\sigma {\omega }^n=(\sigma  {\omega  }^{n-
1}/i)\cdot  i\omega =\Delta ({{\scr H}}_2,n)i\omega $
as given  in
Lemma  (4.10).  By  Proposition (4.11),  the limit  $  \lim_{n\to
\infty }D_n$ exists.  If $\lim_n D_n =\infty $,  then it follows
from ${\xi }_n \ge 1/2 $ that
$$\lim_{n\to \infty }\frac {\delta ({{\scr P}}_2,n)}
{\delta ({{\scr H}}_2,n)} =\lim_{n\to \infty } \frac {\Delta ({\scr
P}_2,n)} {\Delta ({{\scr H}}_2,n)} .\tag 4.14 $$
Hence,  we  may and will assume that $\lim_{n\to \infty }  D_n  <
\infty.$ But then,  from (4.11) and the construction of ${\varphi
}_j$'s, we must have
$$  \frac {\Delta ({{\scr P}}_2,n+1)} {\Delta ({{\scr P}}_2,n)\omega }
=\lambda (n,n+1) \longrightarrow 1 \quad \text{as} \; n\to \infty.
\tag 4.15 $$
On  the  other hand,  as we have mentioned above,  for any  policy
${{\scr P}}_2,$ we have
$$  \left(\matrix \delta({{\scr P}}_2,n)   \\
\Delta({{\scr P}}_2,n) \endmatrix \right) \le \left(\matrix 0 &i \\1  &i
\endmatrix  \right)  \left(\matrix  \delta  ({{\scr P}}_2,n+1)  \\
\Delta ({{\scr P}}_2,n+1) \endmatrix \right). \tag 4.16 $$
In particular,
$$ (i+{\xi }_n) \lambda (n-1,n)\omega
=\biggl (i+\frac {\delta ({{\scr P}}_2,n)} {\Delta ({{\scr P}}_2,n)}\biggr )
\frac  {\Delta
({{\scr P}}_2,n)} {\Delta ({{\scr P}}_2,n-1)}=\frac {\delta ({\scr
P}_2,n)+i\Delta ({{\scr P}}_2,n)} {\Delta ({{\scr P}}_2,n-1)} \ge 1.
$$
Hence
$$ \frac {{\xi }_n} {i\omega } \ge \frac  {1-i\lambda
(n-1,n)\omega } {i\lambda (n-1,n){\omega }^2}. $$  From this  and
(4.15), we get
$$ {\varliminf}_{n\to \infty } {\xi }_n/(i\omega)
\ge (1-i\omega)/(i{\omega }^2) =1. \tag 4.17
$$
On  the other hand,  by (4.16) again,  we have $$  \frac {\delta
({{\scr P}}_2,n)} {\delta ({{\scr H}}_2,n)} \le \frac {i\Delta ({\scr
P}_2,n+1)} {\delta ({{\scr H}}_2,n)} = \frac {\Delta ({{\scr P}}_2,n+1)}
{\Delta ({{\scr H}}_2,n+1)}.  $$
Combining  this  with (4.13) and (4.17), we finally arrive at
$$ \align
\lim_n D_n &\le
   \Big(\lim_n D_n\Big){\varliminf}_n \frac {{\xi}_n}{i\omega}\\
&\le {\varliminf}_n \frac {\delta({{\scr P}}_2,n)}{\delta({{\scr H}}_2,n)}\\
&\le {\varlimsup}_n\frac {\delta ({{\scr P}}_2,n)} {\delta ({{\scr H}}_2,n)}\\
&\le \lim_n \frac {\Delta ({{\scr P}}_2,n+1)} {\Delta({{\scr H}}_2,n+1)}\\
&=\lim_n D_n.
\endalign$$
Therefore,  we claim that
$ \lim_n  {{\xi}_n}/ {(i\omega)} =1$
and hence (4.14) holds. \qed \enddemo

To  conclude this section,  we show that the optimal policy  at
infinity is essentially unique.

\proclaim{(4.18)Corollary}  If ${{\scr P}}_2$ is not  equivalent
to ${{\scr H}}_2$, then we have $\dsize\lim_{n\to \infty } \dfrac {\delta
({{\scr P}}_2,n)} {\delta ({{\scr H}}_2,n)} > 1$. \endproclaim

\demo{Proof} The conclusion follows from Proposition (4.11)
and Theorem (4.12) immediately. \qed \enddemo

\heading 5. Proof of Theorem (1.15). \endheading

We begin this section by introducing a comparison lemma.
Suppose that we are now at the $N$-th step.
Then the policy $ {{\scr P}}_1$ corresponds in
a  natural  way  a basic policy ${{\scr P}}_2$  having  successive testing
numbers  $k_1=2i$, $k_n=2i-1 $ for all $n\ge 2$ and  initial testing interval
$[0,{\Delta }_{N-1}]$. Moreover,
$ \delta ({{\scr P}}_1, N+m-1)=\delta ({{\scr P}}_2,m)$, $m\ge 1$.  Recalling that for the policy  ${{\scr H}}_2$
defined by Lemma (4.10),  we have
$ \delta ({{\scr H}}_2,m)=\sigma {\omega }^m$,  where
$$\sigma = {\Delta }_{N-1}.  \tag 5.1 $$
As  an  application  of  Theorem  (4.12), we obtain
$$ \align \text{\hskip-1em}
\sup_{m\ge 1} F_{m+2}\delta ({{\scr P}}_1,m) &\ge  \lim_{m\to
     \infty  }  F_{m+2}\delta  ({{\scr P}}_1,m)
=\lim_{m\to \infty}F_{N+m+1} \delta ({{\scr P}}_1,N+m-1)\text{\hskip-4em}\tag5.2 \\
\text{\hskip-1em}&=\lim_{m\to \infty }\bigg(\frac
     {F_{N+m+1}}{F_{m+1}} F_{m+1}\delta  ({{\scr P}}_2,m)\bigg)
\!=\!\frac{1}{\omega^N} \lim_{m\to \infty } F_{m+1}\delta ({{\scr P}}_2,m)\text{\hskip-4em}\\
\text{\hskip-1em}&\ge  {\omega }^{-N} \lim_{m\to \infty } F_{m+1}\delta ({{\scr H}}_2,m)
=\sigma  {\omega  }^{-N}   \lim_{m\to   \infty }
F_{m+1}{\omega  }^m.\text{\hskip-4em} \endalign $$
Thus,  we have  proved  the following result:

\proclaim{(5.3)Lemma}  Let $\gamma =F_4{\omega }^3/\lim_{n\to
\infty} F_{n+1} {\omega }^n$. Then, we have $\delta ({{\scr P}}_1) \ge \delta$
provided
$$ \sigma {\omega }^{-N}\ge \gamma.\tag 5.4 $$
\endproclaim

Based on this  lemma, we can now make a complement to Lemma (3.30).

\proclaim{(5.5)Lemma} Let $m\ge 1$ and $ {\Delta }_m=(-1)^{m-1}(G_{m-3}-G_{m-1}{\alpha }_1)/i^m$.
Then  $\delta({{\scr P}}_1)$ $<\delta $ only if
$$(-1)^m{\alpha }_1 <
\{i^m {\omega }^{m+1}\gamma +(-1)^m G_{m-3} \}/G_{m-1}=: B_m. \tag5.6$$
\endproclaim

\demo{Proof}  Applying  Lemma  (5.3) to the case  that  $
N=m+1,$ we obtain
$$\sigma
{\omega }^{-N} = {\Delta }_m/{\omega  }^{m+1}
= (-1)^{m-1} (G_{m-3}-G_{m-1} {\alpha }_1)/(i^m {\omega }^{m+1}). $$
Thus,  $ \sigma {\omega  }^{-N}\ge
\gamma $ is equivalent to
$ (-1)^m{\alpha }_1 \le  \{ i^m {\omega }^{m+1} \gamma +(-1)^m G_{m-3} \}
/G_{m-1}$.
This proves our assertion. \qed \enddemo

Consider the special case that $m=1$. That is
$${\Delta }_1=(1/i -{\alpha }_1)/i. \tag 5.7 $$
Then,  the  condition  (5.6) becomes
$ {\alpha }_1 \ge \frac 1 i -i \gamma \omega^2 >\frac 1 i - \gamma \omega$.
But we have
$$ 1/i - \gamma \omega >  {G_1}/ {G_3}.\tag5.8 $$
The proofs of this and some subsequent elementary inequalities are delayed
to the end of this section for keeping the main line of the proof of
Theorem (1.15).

Because (5.7) holds for any choice of ${\alpha }_1$,  the above
facts  enable us to remove the first two sub-intervals in  (3.32).
Thus,  for  the rest of the proof,  we need only to consider  the
intervals:
$$\biggl (\frac {G_{2n}} {G_{2n+2}},  \frac {F_{2n}}  {F_{2n+2}}\biggr ),
\quad\quad \biggl (\frac {F_{2n}} {F_{2n+2}}, \frac {G_{2n-2}} {G_{2n}}\biggr).
\tag 5.9 $$

Now, we are at the position to complete the proof of Theorem (1.15). Note
that $\mu (0,1)=\Delta ({{\scr P}}_1,1)\dsize /\Delta ({{\scr H}}_1, 1)=i{\Delta }_1
\dsize /\omega $. Thus,
${\alpha }_1>{\omega }^2 \Leftrightarrow \mu (0,1)<1$.
If $\mu (0,1)^{-1}\ge \gamma $, then ${\varphi }_1\ge \gamma $ and there is
nothing to do. We assume that $\mu (0,1)^{-1}<\gamma $. Equivalently,
${\alpha }_1<\frac 1 i - \frac {\omega } {\gamma }$.
But as we will prove later (Lemma (5.20)) that
$$ 1/i -  {\omega }/{\gamma }<  {F_2}/{F_4}.\tag5.10 $$
This means that we do not
need to consider the sub-interval $(F_2/F_4,G_0/G_2)$. Furthermore, by Lemma
(3.22), for
$${\alpha }_1\in \biggl ( \frac {G_{2n}} {G_{2n+2}},\frac {F_{2n}}
{F_{2n+2}} \biggr ),\quad\quad n\ge 1  \tag5.11
$$
or
$$ {\alpha }_1\in \biggl (
\frac {F_{2n}} {F_{2n+2}}, \frac {G_{2n-2}} {G_{2n}} \biggr ),
\quad\quad n\ge 2 \tag 5.12
$$
the formulas of  ${\delta }_m$ and ${\Delta }_m$ given by (3.23) and (3.24)
are available at least for $m\le 2$. In particular,
$\mu (0,1)=
\Delta ({{\scr P}}_1,2) \dsize /\Delta ({{\scr H}}_1,2)={\alpha }_1\dsize / {\omega
}^2$.
Hence, the proof is deduced to consider the case that ${\omega }^2<
{\alpha }_1<\gamma {\omega }^2$. Given such an ${\alpha }_1$, there exists
uniquely an $n_0$ so that one of (5.11) and (5.12) holds. We now discuss these
two cases separately.

(I) Let (5.12) hold for some $n_0\ge 2$. Then by Lemma (3.22) and Lemma
(3.30), we have $F_{2n_0+1}{\delta }_{2n_0-1}\ge \delta $ unless
${\alpha }_1< A_{2n_0-1}$.
We now prove that this is impossible. This follows once we prove
that $A_{2n_0-1}<F_{2n_0}/F_{2n_0+2}$ which contradicts to our assumption.
To do so, noticing that the last inequality is equivalent to
$\delta <F_2F_{2n_0+1}/F_{2n_0+2}$,
by (2.8), we need only to show that
$$\delta < F_2 F_5/F_6. \tag5.13
$$
We will check this in Lemma (5.15).

(II) Let (5.11) hold for some $n_0\ge 1$. Then, we have $F_{2n_0 +2}
{\delta }_{2n_0}\ge \delta $ unless ${\alpha }_1<$$B_{2n_0}$.
But we can
prove that $ B_{2n_0}<G_{2n_0}\dsize /G_{2n_0+2}$.
This again gives us a
contradiction. Actually, the above inequality is equivalent to
$ {\omega }^{2n_0+1}G_{2n_0+2}<(i+1)/\gamma$.
Hence, it suffices to show that
$${\omega }^3G_4<(i+1)/\gamma . \tag5.14
$$
This will be done by Lemma (5.19).
Finally, we conclude our main proof by the following
four lemmas.

\proclaim{(5.15)Lemma} $\delta <F_2 F_5/F_6. $ \endproclaim

\demo{Proof} By (3.4), $\delta =i(2i+3)((i+1)\omega -1)$.
We need only to show that
$ (i+1)\omega <1+\frac {2F_5} {(2i+3)F_6}$.
But the right hand side equals to
$$\align
1+\frac 2 {2i+3}\cdot \frac {2i^2+5i+1} {i(2i^2+7i+4)} &=1+\frac 2 {i(2i+3)}
\biggl (1-\frac {2i+3} {2i^2+7i+4} \biggr )\\
&=1+\frac 2 {i(2i+3)} -\frac 2 {i(2i^2+7i+4)}. \endalign$$
On the other hand,
$$ (i+1)\omega
\!<\! \frac{(i+1)F_6}{F_7}\! =\!\frac {(i+1)(2i^2+7i+4)} {2i^3+9i^2+9i+1}
\!=\!1+\frac {2i+3} {2i^3+9i^2+9i+1} \!<\!1+\frac 1 {i(i+3)}.\text{\hskip-1em}\tag5.16 $$
Thus,  it suffices to show that
$\dfrac 2 {2i^2+7i+4}<\dfrac 2 {2i+3}-
\dfrac 1 {i+3} =\dfrac 3 {2i^2+9i+9}$.
This certainly holds for all $i\ge 2$. \qed \enddemo

\proclaim{(5.17)Lemma } $1/i  -\gamma \omega > G_1/G_3.$\endproclaim

\demo{Proof} Observe that the assertion is equivalent to
$$ \gamma< (i+1)(1+\omega)/(i+2).\tag 5.18$$
By (3.2) and (3.3), this becomes
$\dfrac {i(i+4)(2i+3)((i+1)\omega-1)}{2(i+1)+3i\omega} <
  \dfrac {i+1}{i+2}(1+\omega)$.
That is
$$\align
i(i+2)(i+4)(2i+3)((i+1)\omega-1)&< (i+1)(1+\omega)(2(i+1)+3i\omega)\\
  &=(i+1)(2i+5+2(i+1)\omega).\endalign$$
Or
$$\omega < \frac{(i+1)(2i+5)+i(i+2)(i+4)(2i+3)}{(i+1)\bigl(i(i+2)(i+4)(2i+3)
    \!-\!2(i+1)\bigr)}
         =\frac {2i^4+15i^3+36i^2+31i+5}{2i^5+17i^4+49i^3+56i^2+20i\!-\!2}.$$
Note that the right hand side is greater than $1/(i+1-1/i)$ and $\omega=
(\sqrt{ 1+4/i} -1)/2$. Now, it should be easy to obtain the required assertion.
\qed \enddemo

\proclaim{(5.19)Lemma } ${\omega }^3G_4<(i+1)/\gamma .$ \endproclaim

\demo{Proof} It follows from (5.16) that
$$G_4{\omega }^3=i^2(i^2+3i+1)
{\omega }^3 =i(i^2+3i+1)((i+1)\omega -1)< (i^2+3i+1)/(i+3).$$
Thus, it suffices to show that
$\gamma <\frac {(i+1)(i+3)} {i^2+3i+1} =1+\frac {i+2} {i^2+3i+1}$.
But this follows from  (5.18) and (5.16):
$$ \gamma <\bigl (1- 1/(i+2)\bigr)(1+\omega)
    =1+ \big( (i+1)\omega-1\big) /(i+2) < 1+ 1/(i+1). \qed $$\enddemo

\proclaim{(5.20)Lemma }  $1/i - \omega /\gamma < F_2 /F_4.$   \endproclaim

\demo{Proof } The assertion is the same as follows:
$ \gamma < i(2i+3)\omega /(2i+1)$.
By (5.18), it is enough to show that
$ \dfrac {i+1}{i+2}(1+\omega) \!<\! \dfrac {i(2i+3)}{2i+1}\omega$.
Equivalently,
$ \omega\! >\! \dfrac {2i^2+3i+1}{2i^3+5i^2+3i-1}$.
Note that the right hand side is less than $1/(i+1-1/(2i-1))$, it is now easy
to complete the proof.\qed \enddemo

\enddocument